\documentstyle[11pt]{article}
\title{On the Osculating Spaces of Submanifolds in Euclidean Spaces}
\author{Kostadin Tren\v{c}evski\\
Inst. of Math., Faculty of Natural Sci. and Mathematics,
\\ Ss. Cyril and Methodius Univ. in Skopje, Arhimediva 3, Macedonia\\ 
e-mail: kostadin.trencevski@gmail.com, kostatre@pmf.ukim.mk}
\date{}

\begin{document}
\maketitle

\begin{abstract}
This paper is a continuation of the papers \cite{2,3,4,5,6}. 
In this paper the osculating spaces of arbitrary order 
of a manifold embedded in Euclidean space are considered. 
A better estimation of their dimensions as well as 
the description of its basis are given. 



\end{abstract}

\section{Introduction}
Let us consider an $n$-dimensional manifold $M$ embedded in the 
Euclidean space $R^m$, where $m=n+k$. In \cite{2} it is 
considered the following $k\times k$ matrix 
$$
P^{(1)}_{\alpha \beta } = \sum _{i=1}^{n}\sum _{j=1}^{n}
(Y_{i}\cdot \nabla _{Y_{j}} N_{\alpha })
(Y_{i}\cdot \nabla _{Y_{j}} N_{\beta }) ,\eqno{(1.1)}
$$
where $Y_1,\ldots ,Y_n$ are orthonormal tangent vectors at a chosen
point on $M$ and $N_1,\ldots ,N_k$ are orthonormal vector fields 
which are orthogonal to the tangent space in a neighborhood of the chosen 
point. It is obvious that $P^{(1)}$ does not depend on the choice 
of the orthonormal base. The eigenvalues of $P^{(1)}$ are non-negative 
numbers and its eigenvectors are orthogonal. Let $k_1=rank (P^{(1)})$ 
and let $\lambda^{(1)}_1,\ldots ,\lambda^{(1)}_{k_1}$ be the positive 
eigenvalues and for any eigenvector $(p_1,\ldots ,p_k)$ we consider
the vector $p_1N_1+\ldots +p_kN_k$ and hence we obtain the following 
vectors $N^{(1)}_1,\ldots ,N^{(1)}_{k_1}$ as eigenvectors from the normal 
space. These vectors do not depend on the choice of the base $N_{\alpha}$.
The positive eigenvalues $\lambda^{(1)}_1,\ldots ,\lambda^{(1)}_{k_1}$
are called the {\em first normal curvatures} and the corresponding 
vectors $N^{(1)}_1,\ldots ,N^{(1)}_{k_1}$ are called the {\em first 
normal vectors}. Note that $k_1\le n^2$, because the right side of 
(1.1) is a sum of $n^2$ matrices of rank 1. 
Using the results in \cite{1}, in 
\cite{6} it is proved that the vectors $Y_1,\ldots ,Y_n,N^{(1)}_1,\ldots ,
N^{(1)}_{k_1}$, generate the osculating space (of the first order) 
at the considered point,
which gives the geometrical interpretation of the first normal vectors. 

In \cite{2} are introduced the second normal curvature tensors and 
the second normal vectors as follows. Without loss of generality we 
assume that $N_1=N_1^{(1)},\ldots ,N_{k_1}=N_{k_1}^{(1)}$. Now it is 
$N_i\cdot N_j^{(1)}=0$ for $i>k_1$ and $j\le k_1$. We consider the 
$(k-k_1)\times (k-k_1)$ matrix 
$$P^{(2)}_{\alpha\beta}=\sum_{i=1}^{n+k_1}\sum_{j=1}^n 
(Y_i\cdot \nabla _{Y_j}N_{\alpha})(Y_i\cdot \nabla _{Y_j}N_{\beta}),$$
for $k_1+1\le \alpha ,\beta \le k$, and where we have denoted 
$Y_{n+1}=N_1,\ldots ,Y_{n+k_1}=N_{k_1}$. According to the choice 
of $N^{(1)}_1,\ldots ,N^{(1)}_{k_1}$ we get the following reduced form 
$$P^{(2)}_{\alpha\beta}=\sum_{i=1}^{k_1}\sum_{j=1}^n 
(N^{(1)}_i\cdot \nabla _{Y_j}N_{\alpha})
(N^{(1)}_i\cdot \nabla _{Y_j}N_{\beta}). \eqno{(1.2)}
$$
If $k_2=rank(P^{(2)})=0$ at any point of the submanifold, then the 
manifold locally can be embedded in $n+k_1$-dimensional affine subspace of 
$R^m$. If $k_2>0$, let $(\lambda_1,\ldots ,\lambda_{k-k_1})$ be an 
eigenvector of $P^{(2)}$, then we consider the vector 
$\lambda_1N_{k_1+1}+\ldots +\lambda_{k-k_1}N_k$ as an eigenvector. 
According to this identification, the eigenvectors of $P^{(2)}$ and the 
principal directions do not depend on the choice of the basis $N_{\alpha}$.
The positive eigenvalues $\lambda^{(2)}_1,\ldots ,\lambda^{(2)}_{k_2}$ 
are defined to be the {\em second normal curvatures} and the corresponding 
eigenvectors $N^{(2)}_1,\ldots ,N^{(2)}_{k_2}$ are defined to be the 
{\em second normal vectors}. In \cite{6} it is proved that the vectors 
$Y_1,\ldots ,Y_n,N^{(1)}_1,\ldots ,
N^{(1)}_{k_1},N^{(2)}_1,\ldots ,N^{(2)}_{k_2}$ 
generate the osculating space of the second order at the considered point,
which gives the geometrical interpretation of the second normal vectors. 

Continuing this procedure, 
the normal curvatures and normal vectors of higher degree are introduced 
\cite{2,6}. 
This procedure is finite, since the number $m$ is finite. 
We give only the inductive step for the matrix 
$P^{(l+1)}$. 
Namely, 
$P^{(l+1)}$ is $(k-k_1 -\cdots -k_l )\times (k-k_1 -\cdots -k_l )$ matrix 
given by 
$$
P^{(l+1)}_{\alpha \beta } = \sum _{i=1}^{k_l}\sum _{j=1}^{n}
(N_{i}^{(l)}\cdot \nabla _{Y_{j}}N_{\alpha })
(N_{i}^{(l)}\cdot \nabla _{Y_{j}}N_{\beta }) ,
\eqno{(1.3)}
$$
and $rank P^{(l+1)}\le n \cdot rank P^{(l)}$. By induction of 
$l$ it follows that $rank (P^{(l)})\le n^{l+1}$. 

The vectors 
$$Y_1,\ldots ,Y_n,N^{(1)}_1,\ldots ,
N^{(1)}_{k_1},N^{(2)}_1,\ldots ,N^{(2)}_{k_2},\ldots ,
N^{(l)}_1,\ldots ,N^{(l)}_{k_l}$$ 
generate the osculating space of order $l$ at the considered point. 

\section{Main result} 
We saw in the introduction that 
$$k_r = rank (P^{(r)}_{\alpha \beta}) \le n^{r+1}. $$
In this paper we prove much better inequality 
$$k_r = rank (P^{(r)}_{\alpha \beta}) \le {n+r\choose r+1}.$$

First, let us consider the case $r=1$. In this case the matrix 
$$P^{(1)}_{\alpha\beta} = \sum_{i=1}^n \sum_{j=1}^n 
(Y_i \cdot \nabla_{Y_j}N_{\alpha})(Y_i \cdot \nabla_{Y_j}N_{\beta})$$
is a sum of $n^2$ $(1\le i,j\le n)$ matrices of type $a_ia_j$ of rank 1. 
But since the matrices for the pairs $(i,j)$ and $(j,i)$ are equal, 
we obtain that $P^{(1)}$ is a sum of 
$n+{n\choose 2}={n+1\choose 2}$ matrices of rank equal to 1, and hence 
$k_1\le {n+1\choose 2}$. Indeed, note that 
$$Y_i \cdot \nabla_{Y_j}N_{\alpha}=-
N_{\alpha}\cdot \nabla_{Y_j}Y_i = $$
$$=- N_{\alpha}\cdot (\nabla_{Y_i}Y_j - [Y_i,Y_j]) = 
- N_{\alpha}\cdot \nabla_{Y_i}Y_j = 
Y_j \cdot \nabla_{Y_i}N_{\alpha},$$
where we used that the torsion tensor 
$T(Y_i,Y_j)=\nabla_{Y_i}Y_j-\nabla_{Y_j}Y_i-[Y_i,Y_j]$ is a zero tensor 
and we used that $[Y_i,Y_j]$ is a tangent vector. Hence we have 
equal summands for the pairs $(i,j)$ and $(j,i)$. 

Moreover, according to the definition of $P^{(1)}_{\alpha\beta}$ 
we see that the space generated by the eigenvectors of 
$P^{(1)}$ coincides with the space generated 
by the vectors 
$$\sum_{\alpha =1}^k N_{\alpha}(Y_i\cdot \nabla _{Y_j}N_{\alpha }) = 
-\sum_{\alpha =1}^k N_{\alpha}(N_{\alpha}\cdot \nabla _{Y_j}Y_i) ,$$
i.e. the projection of the vectors $\nabla _{Y_j}Y_i$ on the normal space
of the tangent space. Hence the first osculating space generated by 
$$Y_1,\ldots ,Y_n,N^{(1)}_1,\ldots ,N^{(1)}_{k_1}$$ 
is the same with 
the space generated by the vectors 
$$Y_1,\ldots ,Y_n, \Bigl \{ \nabla _{Y_i}Y_j \Bigr \}\quad (i\le j).$$

Further we will consider the case $r=2$. 

According to the definition 
$$P^{(2)}_{\alpha\beta} = \sum_{i=1}^{k_1}\sum_{j=1}^n
(N^{(1)}_i\cdot \nabla _{Y_j}N_{\alpha})
(N^{(1)}_i\cdot \nabla _{Y_j}N_{\beta})$$
we obtain that the space of eigenvectors of $P^{(2)}$ 
coincides with the space generated by the vectors of type 
$$\sum_{\alpha =k_1+1}^k N_{\alpha}(N_{\alpha}\cdot \nabla_{Y_j}N^{(1)}_i)
\eqno{(2.1)}$$
for $1\le j\le n$ and $1\le i\le k_1$. Note that the 
vectors $N^{(1)}_1,\ldots ,N^{(1)}_{k_1}$ generate the same space as the 
space of projections of $\nabla_{Y_j}Y_i$ on the normal space, i.e. 
the space of vectors 
$$\nabla _{Y_j}Y_i - \sum_{p=1}^n Y_p(Y_p\cdot \nabla _{Y_j}Y_i).$$
By replacing these vectors instead of $N^{(1)}_i$ in (2.1) we obtain that the
space of eigenvectors of $P^{(2)}$ coincides with the space 
generated by the projection of the vectors 
$$\nabla _{Y_j}\nabla _{Y_i}Y_p \quad (1\le i,j,p\le n)$$
over the space orthogonal to 
$Y_1,\ldots ,Y_n,N^{(1)}_1,\ldots ,N^{(1)}_{k_1}$. Similarly as we proved 
that $\nabla _{Y_j}Y_i$ and $\nabla _{Y_i}Y_j$ differ only for a vector 
in the tangent space, we shall prove that the vectors 
$$\nabla _{Y_j}\nabla _{Y_i}Y_p, \; \nabla _{Y_j}\nabla _{Y_p}Y_i,\;
\nabla _{Y_i}\nabla _{Y_j}Y_p,\; \nabla _{Y_i}\nabla _{Y_p}Y_j,\; 
\nabla _{Y_p}\nabla _{Y_j}Y_i,\; \nabla _{Y_p}\nabla _{Y_i}Y_j,$$
differ only for a vector belonging to the space generated by 
$Y_1,\ldots ,Y_n,$ $N^{(1)}_1,\ldots ,N^{(1)}_{k_1}$. 
Indeed, it is sufficient to prove that the differences 
$$\nabla _{Y_j}\nabla _{Y_i}Y_p-\nabla _{Y_i}\nabla _{Y_j}Y_p
\quad \hbox { and } \quad 
\nabla _{Y_j}\nabla _{Y_i}Y_p-\nabla _{Y_j}\nabla _{Y_p}Y_i$$
belong to that space. Namely, using that 
$$R(Y_j,Y_i)Y_p = \nabla_{Y_j}\nabla_{Y_i}Y_p - \nabla_{Y_i}\nabla_{Y_i}Y_p
-\nabla_{[Y_j,Y_i]}Y_p$$
belongs to the tangent space, and $\nabla_{[Y_j,Y_i]}Y_p$ belongs to 
the space generated by 
$Y_1,\ldots ,Y_n, N^{(1)}_1,\ldots ,N^{(1)}_{k_1}$, we obtain that 
$\nabla_{Y_j}\nabla_{Y_i}Y_p - \nabla_{Y_i}\nabla_{Y_i}Y_p$ belongs to the 
space generated by $Y_1,\ldots ,Y_n, N^{(1)}_1,\ldots ,N^{(1)}_{k_1}$. 
The second difference 
$\nabla _{Y_j}\nabla _{Y_i}Y_p-\nabla _{Y_j}\nabla _{Y_p}Y_i$
belongs to the space generated by 
$Y_1,\ldots ,Y_n, N^{(1)}_1,\ldots ,$ $N^{(1)}_{k_1}$ because 
$$\nabla _{Y_j}\nabla _{Y_i}Y_p-\nabla _{Y_j}\nabla _{Y_p}Y_i=
\nabla _{Y_j}([Y_i,Y_p])$$
and $[Y_i,Y_p]$ is a vector from the tangent space. 

So, without loss of generality we can consider those triples $(i,j,p)$ 
such that $1\le i\le j\le p$. But, such triples there are exactly 
${n(n+1)(n+2)\over 3!}$. Hence, $k_2\le {n+2\choose 3}$. Moreover, according
to this discussion we obtain that the second osculating space, i.e. 
generated by 
$Y_1,\ldots ,Y_n,N^{(1)}_1,\ldots ,N^{(1)}_{k_1},
N^{(2)}_1,\ldots ,N^{(2)}_{k_2}$
coincides with the space generated by the vectors 
$$Y_1,\ldots ,Y_n, \Bigl \{ \nabla _{Y_i}Y_j \Bigr \} (i\le j), 
\Bigl \{ \nabla _{Y_i}\nabla _{Y_j}Y_p \Bigr \} (i\le j\le p).$$

We can continue this consideration and using an inductive step analogously 
to the step from $r=1$ to $r=2$, we obtain the main theorem. 

{\bf Theorem.} 
{\em (i) \quad $k_r=dim P^{(r)}\le {n+r\choose r+1}$;

(ii) The space generated by the eigenvectors of $P^{(r)}$ coincides with the 
vector space generated by the projections of the vectors 
$$\nabla_{Y_{i_1}}\nabla_{Y_{i_2}}\cdots \nabla_{Y_{i_{r-1}}} Y_{i_r}
\qquad (i_1\le i_2\le \cdots \le i_r)$$
to the orthogonal complement of the space generated by the vectors 
$$Y_1,\ldots ,Y_n,N^{(1)}_1,\ldots ,N^{(1)}_{k_1},\ldots ,
N^{(r-1)}_1,\ldots ,N^{(r-1)}_{k_{r-1}};$$

(iii) The r-dimensional osculating space, i.e. the space generated by 
$$Y_1,\ldots ,Y_n,N^{(1)}_1,\ldots ,N^{(1)}_{k_1},\ldots ,
N^{(r)}_1,\ldots ,N^{(r)}_{k_{r}}$$
coincides with the space generated by 
$$Y_1,\ldots ,Y_n, \Bigl \{ \nabla _{Y_i}Y_j \Bigr \} (i\le j), 
\Bigl \{ \nabla _{Y_i}\nabla _{Y_j}Y_p \Bigr \} (i\le j\le p),\ldots ,$$
$$\Bigl \{ \nabla _{Y_{i_1}}\nabla _{Y_{i_2}}\cdots 
\nabla _{Y_{i_{r-1}}}Y_{i_r} \Bigr \} (i_1\le i_2\le \ldots \le i_r),$$
where $Y_1,\ldots ,Y_n$ is an arbitrary orthonormal basis of the tangent 
space. }

At the end we will prove that there does not exist a better 
estimation of $k_r$ than $k_r=dim P^{(r)}\le {n+r\choose r+1}$.

Let $r$ be a given positive integer. We choose $m$ sufficiently large 
number such that $m\ge n+{n+1\choose 2}+{n+2\choose 3}+\cdots +
{n+r\choose r+1}$. Then we define an $n$-dimensional surface in $R^m$
parameterized by $u_1,\ldots ,u_n$ by 
$$f(u_1,\ldots ,u_n)=(u_1,\ldots ,u_n,
a_{11}u_1^2,a_{12}u_1u_2,\ldots ,a_{nn}u_n^2,a_{111}u_1^3,$$
$$a_{112}u_1^2u_2,\ldots ,a_{nnn}u_n^3,\ldots ,
a_{11...1}u_1^r,a_{11..12}u_1^{r-1}u_2,\ldots ,a_{nn...n}u_n^r,0,\ldots ,0),
$$
where all of the coefficients 
$a_{i_1i_2}$, $a_{i_1i_2i_3}$,...,$a_{i_1i_2\ldots i_r}$
$(i_1\le i_2\le \ldots \le i_r)$ are nonzero coefficients. In this case (at the coordinate origin) we have 
$$k_1= {n+1\choose 2}, \quad k_2={n+2\choose 3},\quad \ldots , 
\quad k_r= {n+r\choose r+1}.$$

\end{document}